\documentclass[12pt]{article}
\usepackage{amsthm,amsfonts}
\usepackage{fullpage}

\newcommand{\binary}{\lbrace 0,1 \rbrace}

\newtheorem{theorem}{Theorem}

\newtheorem{lemma}[theorem]{Lemma}
\newtheorem{proposition}[theorem]{Proposition}

\title{Binary words containing infinitely many overlaps}
\author{James Currie \\
Department of Mathematics \\
University of Winnipeg \\
Winnipeg, Manitoba R3B 2E9 (Canada) \\
{\tt j.currie@uwinnipeg.ca} \bigskip \\
Narad Rampersad, Jeffrey Shallit \\
School of Computer Science \\
University of Waterloo \\
Waterloo, Ontario N2L 3G1 (Canada) \\
{\tt nrampersad@math.uwaterloo.ca} \\
{\tt shallit@graceland.math.uwaterloo.ca} \\
}

\begin{document}
\date{\today}
\maketitle

\begin{abstract}
We characterize the squares occurring in infinite overlap-free
binary words and construct various $\alpha$~power-free binary words
containing infinitely many overlaps.
\end{abstract}

\section{Introduction}
If $\alpha$ is a rational number, a word $w$ is an $\alpha$~\emph{power}
if there exists words $x$ and $x'$, with $x'$ a prefix of $x$, such that
$w = x^nx'$ and $\alpha = n + |x'|/|x|$.  We refer to $|x|$ as a
\emph{period} of $w$.  An $\alpha^+$~\emph{power}
is a word that is a $\beta$~power for some $\beta > \alpha$.
A word is $\alpha$~\emph{power-free} (resp. $\alpha^+$~\emph{power-free})
if none of its subwords is an $\alpha$~power (resp. $\alpha^+$~power).
A $2$~power is called a \emph{square}; a $2^+$~power is called an
\emph{overlap}.

Thue \cite{Thu12} constructed an infinite overlap-free binary word;
however, Dekking \cite{Dek76} showed that any such infinite word
must contain arbitrarily large squares.  Shelton and Soni \cite{SS85}
characterized the overlap-free squares, but it is not hard to show that
there are some overlap-free squares, such as $00110011$, that cannot
occur in an infinite overlap-free binary word.  In this paper, we
characterize those overlap-free squares that do occur in infinite
overlap-free binary words.

Shur \cite{Shu00} considered the bi-infinite overlap-free and
$7/3$~power-free binary words and showed that these classes of words were
identical.  There have been several subsequent papers
\cite{AC04,KS03,KKT99,Ram05} that have shown various similarities
between the classes of overlap-free binary words and
$7/3$~power-free binary words.  Here we contrast the two classes of
words by showing that there exist one-sided infinite $7/3$~power-free binary
words containing infinitely many overlaps.  More generally, we show that
for any real number $\alpha > 2$ there exists a real number $\beta$
arbitrarily close to $\alpha$ such that there exists an infinite
$\beta^+$~power-free binary word containing infinitely many $\beta$~powers.

All binary words considered in the sequel will be over the alphabet
$\binary$.  We therefore use the notation $\overline{w}$ to denote
the \emph{binary complement} of $w$; that is, the word obtained from
$w$ by replacing $0$ with $1$ and $1$ with $0$.

\section{Properties of the Thue-Morse morphism}
In this section we present some useful properties of the
\emph{Thue-Morse morphism}; \emph{i.e.}, the morphism $\mu$
defined by $\mu(0) = 01$ and $\mu(1) = 10$.
It is well-known \cite{MH44,Thu12} that the \emph{Thue-Morse word}
\[\mathbf{t} = \mu^\omega(0) = 0110100110010110\cdots\]
is overlap-free.

The following property of $\mu$ is easy to verify.

\begin{lemma}
\label{tmmorph}
Let $x$ and $y$ be binary words.  Then $x$ is a prefix (resp. suffix)
of $y$ if and only if $\mu(x)$ is a prefix (resp. suffix) of $\mu(y)$.
\end{lemma}

Shur \cite{Shu00} proved the following useful theorem.

\begin{theorem}[Shur]
\label{shur}
Let $w$ be a binary word and let $\alpha > 2$ be a real number.  Then $w$
is $\alpha$~power-free if and only if $\mu(w)$ is $\alpha$~power-free.
\end{theorem}

The following sharper version of one direction of this theorem (implicit in
\cite{KS03}) is also useful.

\begin{theorem}\label{stronger}
Suppose $\mu(w)$ contains a subword $u$ of period $p$, with
$|u|/p>2$.  Then $w$ contains a subword $v$ of length
$\lceil|u|/2\rceil$ and period $p/2$.
\end{theorem}

Karhum\"aki and Shallit \cite{KS03} gave the following generalization of
the factorization theorem of Restivo and Salemi \cite{RS84}.  The extension
to infinite words is clear.

\begin{theorem}[Karhum\"aki and Shallit]
\label{fact}
Let $x\in\binary^*$ be $\alpha$~power-free, $2 < \alpha \leq 7/3$.
Then there exist $u,v\in\{\epsilon,0,1,00,11\}$ and an
$\alpha$~power-free $y\in\binary^*$ such that $x=u\mu(y)v$.
\end{theorem}

\section{Overlap-free squares}
Let \[A = \{00,11,010010,101101\}\] and let
\[\mathcal{A} = \bigcup_{k \geq 0}\mu^k(A).\]
Pansiot \cite{Pan81} and Brlek \cite{Brl89} gave the following
characterization of the squares in $\mathbf{t}$.

\begin{theorem}[Pansiot; Brlek]
\label{pansiot}
The set of squares in $\mathbf{t}$ is exactly the set $\mathcal{A}$.
\end{theorem}

We can use this result to prove the following.

\begin{proposition}
\label{square}
For any position $i$, there is at most one square in $\mathbf{t}$
beginning at position $i$.
\end{proposition}

\begin{proof}
Suppose to the contrary that there exist distinct squares $x$ and $y$
that begin at position $i$.  Without loss of generality, suppose that
$x$ and $y$ begin with $0$.  Then by Theorem~\ref{pansiot},
$x = \mu^p(u)$ and $y = \mu^q(v)$, for some $p,q$ and
$u,v \in \{00,010010\}$.  Suppose $p \leq q$ and let $w = \mu^{q-p}(v)$.
By Lemma~\ref{tmmorph}, either $u$ is a proper prefix of $w$ or $w$ is a
proper prefix of $u$, neither of which is possible for any choice of
$u,v \in \{00,010010\}$.
\end{proof}

The set $\mathcal{A}$ does not contain all possible overlap-free squares.
Shelton and Soni \cite{SS85} characterized the overlap-free squares
(the result is also attributed to Thue in \cite{Ber92}).

\begin{theorem}[Shelton and Soni]
\label{conj}
The overlap-free binary squares are the conjugates of the
words in $\mathcal{A}$.
\end{theorem}

Some overlap-free squares cannot occur in any infinite overlap-free
binary word, as the following lemma shows.

\begin{lemma}
\label{extend}
Let $x = \mu^k(z)$ for some $k \geq 0$ and $z \in \{011011,100100\}$.
Then $xa$ contains an overlap for all $a\in\binary$.
\end{lemma}

\begin{proof}
It is easy to see that $x = uvvuvv$ for some $u,v\in\binary^*$,
where $u$ and $v$ begin with different letters.  Thus
one of $uvvuvva$ or $vva$ is an overlap.
\end{proof}

We can characterize the squares that can occur in an infinite
overlap-free binary word.  Let \[B = \{001001,110110\}\] and let
\[\mathcal{B} = \bigcup_{k \geq 0}\mu^k(B).\]

\begin{theorem}
The set of squares that can occur in an infinite overlap-free binary word
is $\mathcal{A} \cup \mathcal{B}$. Furthermore,
if $\mathbf{w}$ is an infinite overlap-free binary word containing
a subword $x\in\mathcal{B}$, then $\mathbf{w}$ begins with $x$
and there are no other occurrences of $x$ in $\mathbf{w}$.
\end{theorem}

\begin{proof}
Let $\mathbf{w}$ be an infinite overlap-free binary word beginning
with a square $yy \not\in \mathcal{A} \cup \mathcal{B}$.  Suppose further that
$yy$ is a smallest such square that can be extended to an infinite
overlap-free word.  If $|y| \leq 3$, then
$yy \not\in \mathcal{A} \cup \mathcal{B}$ is one of $011011$
or $100100$, neither of which can be extended to an infinite
overlap-free word by Lemma~\ref{extend}.

We assume then that $|y| > 3$.  Since, by Theorem~\ref{conj}, $yy$ is a
conjugate of a word in $\mathcal{A}$, we have two cases.

Case 1: $yy = \mu(zz)$ for some $z\in\binary^*$.  By Theorem~\ref{fact},
$\mathbf{w} = \mu(zz\mathbf{w}')$ for some infinite $\mathbf{w}'$, where
$zz\mathbf{w}'$ is overlap-free.  Thus $zz$ is a smaller square
not in $\mathcal{A} \cup \mathcal{B}$ that can be extended to an
infinite overlap-free word, contrary to our assumption.

Case 2: $yy = a\mu(zz')\overline{a}$ for some $a\in\binary$ and
$z,z'\in\binary^*$.  By Theorem~\ref{fact}, $yy$ is followed by
$a$ in $\mathbf{w}$, and so $yya$ is an overlap, contrary to our assumption.

Since both cases lead to a contradiction, our assumption that
$yy \not\in \mathcal{A} \cup \mathcal{B}$ must be false.

To see that each word in $\mathcal{A} \cup \mathcal{B}$ does occur
in some infinite overlap-free binary word, note that Allouche, Currie, and
Shallit \cite{ACS98} have shown that the word
$\mathbf{s} = 001001\mathbf{\overline{t}}$ is overlap-free.  Now consider
the words $\mu^k(\mathbf{s})$ and $\mu^k(\mathbf{\overline{s}})$, which are
overlap-free for all $k \geq 0$.

Finally, to see that any occurrence of $x \in \mathcal{B}$
in $\mathbf{w}$ must occur at the beginning of $\mathbf{w}$, we note that
by an argument similar to that used in Lemma~\ref{extend}, $ax$ contains
an overlap for all $a\in\binary$, and so $x$ occurs at the beginning of
$\mathbf{w}$.
\end{proof}

\section{Words containing infinitely many overlaps}
In this section we construct various infinite $\alpha$~power-free binary
words containing infinitely many overlaps.  We begin by considering the
infinite $7/3$~power-free binary words.

\begin{proposition}
For all $p \geq 1$, an infinite $7/3$~power-free word
contains only finitely many occurrences of overlaps with period $p$.
\end{proposition}

\begin{proof}
Let $\mathbf{x}$ be an infinite $7/3$~power-free word containing
infinitely many overlaps with period $p$.  Let $k \geq 0$ be the smallest
integer satisfying  $p \leq 3 \cdot 2^k$.  Suppose $\mathbf{x}$ contains
an overlap $w$ with period $p$ starting in a position $\geq 2^{k+1}$.
Then by  Theorem~\ref{fact}, we can write
\[
\mathbf{x} = u_1 \mu(u_2) \cdots \mu^{k-1}(u_k) \mu^k(\mathbf{y}),
\]
where each $u_i \in \{\epsilon,0,1,00,11\}$.  The overlap $w$ occurs
as a subword of $\mu^k(\mathbf{y})$.  By Lemma~\ref{stronger},
$\mathbf{y}$ contains an overlap with period $p / 2^k \leq 3$.
But any overlap with period $\leq 3$ contains a $7/3$~power.  Thus,
$\mathbf{x}$ contains a $7/3$~power, a contradiction.
\end{proof}

The following theorem provides a striking contrast to Shur's result
\cite{Shu00} that the bi-infinite $7/3$~power-free words are overlap-free.

\begin{theorem}
\label{infinite}
There exists a $7/3$~power-free binary word containing
infinitely many overlaps.
\end{theorem}

\begin{proof}
We define the following sequence of words: $A_0 = 00$ and
$A_{n+1} = 0\mu^2(A_n)$, $n \geq 0$.  The first few terms in this sequence
are
\begin{eqnarray*}
A_0 & = & 00 \\
A_1 & = & 001100110 \\
A_2 & = & 0011001101001100101100110100110010110 \\
& \vdots &
\end{eqnarray*}

We first show that in the limit as $n \rightarrow \infty$,
this sequence converges to an infinite word $\mathbf{a}$.
It suffices to show that for all $n$, $A_n$ is
a prefix of $A_{n+1}$.  We proceed by induction on $n$.  Certainly,
$A_0 = 00$ is a prefix of $A_1 = 0\mu^2(00) = 001100110$.
Now $A_n = 0\mu^2(A_{n-1})$, $A_{n+1} = 0\mu^2(A_n)$, and by
induction, $A_{n-1}$ is a prefix of $A_{n}$.  Applying Lemma~\ref{tmmorph},
we see that $A_n$ is a prefix of $A_{n+1}$, as required.

Note that for all $n$, $A_{n+1}$ contains $\mu^{2n}(A_1)$ as a subword.
Since $A_1$ is an overlap with period $4$, $\mu^{2n}(A_1)$ contains
$2^{2n}$ overlaps with period $2^{2n+2}$.  Thus, $\mathbf{a}$ contains
infinitely many overlaps.

We must show that $\mathbf{a}$ does not contain a $7/3$~power.  It suffices
to show that $A_n$ does not contain a $7/3$~power for all $n \geq 0$.
Again, we proceed by induction on $n$.  Clearly, $A_0 = 00$ does not
contain a $7/3$~power.  Consider $A_{n+1} = 0\mu^2(A_n)$.  By induction,
$A_n$ is $7/3$~power-free, and by Theorem~\ref{shur}, so is $\mu^2(A_n)$.
Thus, if $A_{n+1}$ contains a $7/3$~power, such a $7/3$~power must
occur as a prefix of $A_{n+1}$.  Note that $A_{n+1}$ begins with
$00110011$.  The word $00110011$ cannot occur anywhere else
in $A_{n+1}$, as that would imply that $A_{n+1}$ contained a cube $000$ or
$111$, or the $5/2$~power $1001100110$.  If $A_{n+1}$ were to begin with
a $7/3$~power with period $\geq 8$, it would contain two occurrences of
$00110011$, contradicting our earlier observation.  We conclude that
the period of any such $7/3$~power is less than $8$.  Checking that no such
$7/3$~power exists is now a finite check and is left to the reader.
\end{proof}

In fact, we can prove the following stronger statement.

\begin{theorem}
\label{uncount}
There exist uncountably many $7/3$~power-free binary words containing
infinitely many overlaps.
\end{theorem}

\begin{proof}
For a finite binary sequence $b$, we define an operator $g_b$ on
binary words recursively by
\begin{eqnarray*}
g_\epsilon(w)&=&w\\
g_{0b}(w)&=&\mu^2(g_b(w))\\
g_{1b}(w)&=&0\mu^2(g_b(w)).
\end{eqnarray*}
Note that $g_b(0)$ always starts with a 0, so that for any finite
binary words $p$ and $b$, $g_p(0)$ is always a prefix of
$g_{pb}(0)$. Since $g_0(0)$ is not a prefix of $g_1(0)$, $g_{p0}(0)$
is not a prefix of $g_{p1}(0)$ for any $p$, so that distinct $b$
give distinct words. Given an infinite binary sequence
$\mathbf{b}=b_1b_2b_3\cdots$ where the $b_i\in\{0,1\}$, define an infinite
binary sequence $w_{\bf b}$ to be the limit of
$$g_\epsilon(00),g_{b_1}(00), g_{b_1b_2}(00),
g_{b_1b_2b_3}(00),\ldots$$

By an earlier argument, each $w_{\bf b}$ is $7/3$~power-free.
Since $g_1(00) = 001100110$ is an overlap, $g_{b1}(00) = g_b(001100110)$
ends with an overlap for any finite word $b$.  Thus, each 1 in
$\mathbf{b}$ introduces an overlap in $w_{\bf b}$. Since uncountably many
binary sequences contain infinitely many 1's, uncountably many of the
$w_{\bf b}$ are $7/3$~power-free words containing infinitely many
overlaps.
\end{proof}

Next, we show that the sequence $\mathbf{a}$ constructed in the proof
of Theorem~\ref{infinite} is an automatic sequence (in the sense of
\cite{AS03}).

\begin{proposition}
The sequence $\mathbf{a}$ is $4$-automatic.
\end{proposition}

\begin{proof}
We show that $\mathbf{a} = g(h^\omega(0))$, where $h$ and $g$ are the
morphisms defined by

\begin{minipage}{2.5in}
\begin{eqnarray*}
h(0) & = & 0134 \\
h(1) & = & 2134 \\
h(2) & = & 3234 \\
h(3) & = & 2321 \\
h(4) & = & 3421
\end{eqnarray*}
\end{minipage}
\hfill\parbox{0.5in}{and}\hfill
\begin{minipage}{2.5in}
\begin{eqnarray*}
g(0) & = & 0 \\
g(1) & = & 0 \\
g(2) & = & 0 \\
g(3) & = & 1 \\
g(4) & = & 1.
\end{eqnarray*}
\end{minipage}
\bigskip

We make some observations concerning 2-letter subwords: The
sequence $h^\omega(0)$ clearly does not contain any of the words
11, 14, 22, 24, 31, 33, 41 or 44. In fact, neither 12 nor 43
appears as a subword either: Words 12 and 43 do not appear
internally in $h(i)$, $0\le i\le 4$; therefore, if $43$ appears in
$h^n(0)$, it must `cross the boundary' in one of $h(12)$, $h(14)$,
$h(22)$ or $h(24)$. Since 14, 22 and 24 do not appear in
$h^\omega(0)$, word 43 can only appear in $h^n(0)$ as a descendant
of a subword 12 in $h^{n-1}(0)$. However, the situation is
symmetrical; word 12 can only appear in $h^n(0)$ as a descendant
of a subword 43 in $h^{n-1}(0)$. By induction, neither 43 nor 12
ever appears.

The point of the previous paragraph is that
\begin{eqnarray*}
h(0)\mbox{ always occurs in the context }h(0)2\\
h(1)\mbox{ always occurs in the context }h(1)2\\
h(2)\mbox{ always occurs in the context }h(2)2\\
h(3)\mbox{ always occurs in the context }h(3)3\\
h(4)\mbox{ always occurs in the context }h(4)3
\end{eqnarray*}

The word $h^\omega(0)$ can thus be parsed in terms of a new
morphism $f$:
\begin{eqnarray*}
f(0) & = & 1342 \\
f(1) & = & 1342 \\
f(2) & = & 2342 \\
f(3) & = & 3213 \\
f(4) & = & 4213.
\end{eqnarray*}

The parsing in terms of $f$ works as follows: If we write
$h^\omega(0) = 0w$, then $w=f(0w)$. It is useful to rewrite this
relation in terms of the finite words $h^n(0)$. For non-negative
integer $n$ let $x_n$ be the unique letter such that $h^n(0)x_n$
is a prefix of $h^\omega(0)$. Thus $x_0=1$, $x_1=2$, etc. We then
have
\begin{equation}\label{recursion for h}
h^n(0)x_n = 0f(h^{n-1}(0)), \quad n\ge 1.
\end{equation}

Since for all $a \in \{0,1,2,3,4\}$, $g(f(a)) = \mu^2(g(a))$, we have
$g(f(u))=\mu^2(g(u))$ for all words $u$. Therefore, applying $g$ to
(\ref{recursion for h})
\begin{eqnarray*}
g(h^n(0)x_n) &=& g(0f(h^{n-1}(0)))\\
&=& g(0)g(f(h^{n-1}(0)))\\
&=& 0\mu^2(g(h^{n-1}(0))), \quad n\ge 1.
\end{eqnarray*}

From this relation we show by induction that $A_n$ is the prefix
of $g(h^{n+1}(0))$ of length $(4^{n+1}+3\cdot 4^n -1)/3$.
Certainly, $A_0 = 00$ is the prefix of length $2$ of
$g(h(0)) = 0011$.  Consider $A_n = 0\mu^2(A_{n-1})$.
We can assume inductively that $A_{n-1}$ is the prefix of
$g(h^n(0))$ of length $(4^n+3\cdot 4^{n-1} -1)/3$.
Writing $g(h^n(0)) = A_{n-1}z$ for some $z$, we have
\begin{eqnarray*}
g(h^{n+1}(0)x_{n+1}) & = & 0\mu^2(g(h^n(0))) \\
& = & 0\mu^2(A_{n-1}z) \\
& = & A_n \mu^2(z),
\end{eqnarray*}
for some $x_{n+1}$, whence $A_n$ is a prefix of $g(h^{n+1}(0))$.
Since $|A_n| = 4|A_{n-1}| + 1$, we have $|A_n| = (4^{n+1}+3\cdot 4^n -1)/3$,
as required.
\end{proof}

The result of Theorem~\ref{infinite} can be strengthened even further.
\begin{theorem}
For every real number $\alpha > 2$ there exists a real number
$\beta$ arbitrarily close to $\alpha$, such that there is an
infinite $\beta^+$~power-free binary word containing infinitely
many $\beta$~powers.
\end{theorem}

\begin{proof}
Let $s\ge 3$ be a positive integer, and let $r=\lfloor \alpha +
1\rfloor$. Let $t$ be the largest positive integer such that
$r-t/2^s>\alpha$, and such that the word obtained by removing a
prefix of length $t$ from $\mu^s(0)$ begins with 00. Let
$\beta=r-t/2^s.$ Since $\alpha\ge r-1$, we have $t<2^s$. Also,
$\mu^3(0)=01101001$ and $\mu^3(1)=10010110$ are of length 8, and
both contain 00 as a subword; it follows that $|\alpha-\beta|\le
8/2^s$, so that by choosing large enough $s$, $\beta$ can be made
arbitrarily close to $\alpha$.

We construct sequences of words $A_n$, $B_n$ and $C_n$. Define
$C_0=00$. For each $n
\ge 0$:

\begin{enumerate}
\item Let $A_n = 0^{r-2}C_n$.
\item Let $B_n = \mu^{s}(A_n)$.
\item Remove the first $t$ letters from $B_n$ to obtain a new word
$C_{n+1}$ beginning with $00$.
\end{enumerate}

Since each $A_n$ begins with the $r$ power $0^r$, each
$B_n=\mu^s(A_n)$ begins with an $r$ power of period $2^s$.
Removing the first $t$ letters ensures that $C_{n+1}$ commences
with an $(r2^s-t)/2^s$ power, viz., a $\beta$ power. The limit of
the $C_n$ gives the desired infinite word. Let us check that this
limit exists:

Let $w$ be the word consisting of the first $t$ letters of
$\mu^s(0)$. Since all the $A_n$ commence with 0 by construction,
all the $B_n$ commence with $\mu^s(0)$, and hence with $w$. This
means that $B_n=wC_{n+1}$ for each $n$.

We show that $A_n$ is always a prefix of $A_{n+1}$ by induction.
Certainly $A_0$ is a prefix of $A_1$. Assume that $A_{n-1}$ is a
prefix of $A_n$. Since $A_n = 0^{r-2}C_n$ and $A_{n+1} =
0^{r-2}C_{n+1}$, $A_n$ is a prefix of $A_{n+1}$ if $C_n$ is a
prefix of $C_{n+1}$. Since $B_{n-1} = wC_n$ and $B_n = wC_{n+1}$,
$C_n$ is a prefix of $C_{n+1}$ if $B_{n-1}$ is a prefix of $B_n$.
By Lemma~\ref{tmmorph}, $B_{n-1}$ is a prefix of $B_n$ if
$A_{n-1}$ is a prefix of $A_n$, which is our inductive assumption.
We conclude that $A_n$ is a prefix of $A_{n+1}$.

It follows that $C_n$ is a prefix of $C_{n+1}$ for $n \geq 0$, so
that the limit of the $C_n$ exists. It will thus suffice to prove
the following claim: \vspace{.1in}

\noindent{\bf Claim:} The $A_n$, $B_n$ and $C_n$ satisfy the
following:
\begin{enumerate}
\item{The word $C_n$ contains no $\beta^+$ powers.}
\item{The only $\beta^+$ power in $A_n$ is $0^r$.}
\item{Any $\beta^+$ powers in $B_n$ appear only in the prefix
$\mu^s(0^r)$.}

\end{enumerate}

Certainly $C_0$ contains no $\beta^+$ powers, and since
$\beta>r-1$, the only $\beta^+$ power in $A_0$ is $0^r$. Suppose
then that the claim holds for $A_n$ and $C_n$.

Now suppose that $B_n=\mu^s(0^{r-2})\mu^s(C_n)$ contains a
$\beta^+$ power $u$ with period $p$. Since $C_n$ contains no
$\beta^+$ powers, Theorem~\ref{shur} ensures that $\mu^s(C_n)$
contains no $\beta^+$ powers. We can therefore write $B_n=xuy$
where $|x|<|\mu^s(0^{r-2})|$. In other words, $u$ overlaps
$\mu^s(0^{r-2})$ from the right. By Theorem~\ref{stronger}, the
preimage of $B_n$ under $\mu$, i.e., $\mu^{s-1}(A_n)$, contains a
$\beta^+$ power of length at least $|u|/2$ and period $p/2$. In
fact, iterating this argument, $A_n$ contains a $\beta^+$ power of
period $p/2^s$ of length at least $|u|/2^s$. Since the only
$\beta^+$ power in $A_n$ is $0^r$, with period 1, we see that
$p/2^s=1$, whence $p=2^s$ and $|u|\le r2^s$.

Recall that $B_n$ has a prefix $\mu^s(0^{r})$ which also has
period $2^s$, and that this prefix is overlapped by $u$. It
follows that all of $xu$ is a $\beta^+$ power with period $p=2^s$.
However, as just argued, this means that $|xu|\le
r2^s=|\mu^s(0^{r})|$, so that $u$ is contained in $\mu^s(0^r)$ and
part 3 of our claim holds for $B_n$. We now show that parts 1 and
2 hold for $C_{n+1}$ and $A_{n+1}$ respectively, and the truth of
our claim will follow by induction.

Part 1 follows immediately from part 3.

Now suppose that $A_{n+1}$ contains a $\beta^+$ power $u$. Recall
that $A_{n+1}=0^{r-2}C_{n+1}$, and $C_{n+1}$ begins with $00$, but
contains no $\beta^+$ powers. It follows that $u$ is not a subword
of $C_{n+1}$. Therefore, $000$ must be a prefix of $u$. If $u=0^q$
for some integer $q$, then $q\le r$ by the construction of
$A_{n+1}$, and
$$r\ge q>\beta>\alpha>r-1.$$ This implies that $q=r$, and $u=0^r$, as
claimed. If we cannot write $u=0^q$, then $|u|_1\ge 1$. Because
$u$ is a $2^+$ power, 000 must appear twice in $u$ with a 1 lying
somewhere between the two appearances. This implies that $000$ is
a subword of $C_{n+1}$, and hence of $B_n=\mu^s(A_n)$. However, no
word of the form $\mu(w)$ contains $000$. This is a contradiction.
\end{proof}

We conclude by presenting the following open problem.

\begin{quote}
Does there exist a characterization (in the sense of \cite{Ber94,Fif80})
of the infinite $7/3$~power-free binary words?
\end{quote}

\end{document}